\documentclass[12pt]{article}
 
\usepackage{amstext,amssymb,amsmath,stmaryrd,amsthm}
\usepackage{geometry} 
\usepackage{color}
\usepackage{hyperref}
\usepackage{makeidx}
\makeindex

\numberwithin{equation}{section}

\newtheorem{theorem}{Theorem}[section]
\newtheorem{corollary}[theorem]{Corollary}

\newtheorem{proposition}[theorem]{Proposition}

\newtheorem{Remark}[theorem]{Remark}

\begin{document}

\title{A note on the Modified Log-Sobolev inequality.}

\author{Ioannis Papageorgiou\thanks{Institut de Math\'ematiques de Toulouse, Universit\'e Paul Sabatier, 31062 Toulouse cedex 09, France. Email: ioannis.papageorgiou@math.univ-toulouse.fr}}


 \date{}

\maketitle

\begin{abstract}
A criterion is presented for the Modified Logarithmic Sobolev inequality on metric measure spaces. The criterion based on U-bound inequalities introduced by Hebisch and Zegarlinski  allows to show the inequality for measures that go beyond log-concavity.  

~


\end{abstract}

 \section{Introduction.}
  A lot of attention has been focused on inequalities that interpolate between the  Log-Sobolev inequality $$Ent_{\mu}(\vert f\vert^r)\leq C\mu \left\vert \nabla f\right \vert ^r  \   \   \   \    \   \  \   \    (LS_r)$$
  (see \cite{G}, \cite{Ba}, \cite{B-Z}) for $r=2$ and the Spectral Gap inequality, that is
$$\mu \vert f-\mu f\vert^2\leq C'\mu \left\vert \nabla f\right\vert^2 \   \   \   \   \   (SG)$$
  where for the measure  $\mu$ on $\mathbb{R}^n$, the entropy$$Ent_\mu (f):=\mu\left(f\log\frac{f}{\mu f}\right)$$ and $\left\vert \nabla f\right \vert $ is the Euclidean length of the gradient $\nabla f$ of the function $f:\mathbb{R}^n \rightarrow\mathbb{R}$.  For a detailed account of these developments one can look on \cite{G-G-M1} and \cite{B-R2}.  A first example of an  inequality interpolating between the Log-Sobolev and the Spectral Gap was introduced  by \cite{Be} and then studied by \cite{L-O}  and \cite{B-R1}. Modified Log-Sobolev inequalities have properties similar with the Log-Sobolev inequality, since they are both preserved under bounded perturbations and  product measures (see \cite{A-B-C}, \cite{B-Z}), while they both imply concentration of measure properties (see \cite{L}, \cite{G-G-M1}, \cite{B-C-R}, \cite{B-R2}, \cite{P}).

In this paper we are interested in the Modified Log-Sobolev inequality  introduced by \cite{G-G-M1}, (see also \cite{G-G-M2} and \cite{B-R2}).

 Given   a $q\in (2,+\infty)$, for  $x \in \mathbb{R}$ we can define the \textit{modification} $H_q$ of $\vert x\vert ^q:\mathbb{R}\rightarrow\mathbb{R}$ to be
\begin{align*}H_{q}(x)=\begin{cases}x^{2} & \text{if \ } \vert x\vert \leq 1 \\
\vert x \vert^{q} & \text{if \ }\vert x\vert \geq\ 1 \end{cases}
\end{align*}
The definition of the  MLS$(H_q)$ inequality follows.

\noindent
\textit{The Modified  Log-Sobolev  Inequality.} 

\noindent
We say that the measure $\mu$ satisfies the Modified Log-Sobolev Inequality if there exists a constant $C_{MLS}$   such that for any function $f\in C^\infty$  the following holds
$$\mu\vert f\vert^2 log\frac{\vert f\vert^2}{\mu\vert f\vert^2}\leq C_{MLS}
\int   H_q\left(\frac{ \vert \nabla f\vert}{f}\right)f^2d\mu  \   \   \   \   \   \   \ (MLS(H_q))$$
  for some positive  constant  $C_{MLS}$.

Concerning measures that satisfy an MLS($H_q$) inequality, in \cite{G-G-M1} it was shown that the (MLS$(H_q)$) inequality for the one dimensional measure on $\mathbb{R}$ is satisfied by the measure 
$$\mu_p=\frac{e^{-\vert x \vert^p}}{Z_p}dx$$
for $p$ conjugate of $q$, i.e. $\frac{1}{p}+\frac{1}{q}=1$, and $Z_p$ the normalization constant. Furthermore, Barthe and Roberto in \cite{B-R2}  showed that a measure $\frac{e^{-V(x)}}{\int e^{-V(x)}dx}dx$ such that 
 $$\liminf_{\vert x \vert\rightarrow \infty}sign(x)V'(x)>0 \text{ \   and  \  }\lim_{\vert x \vert\rightarrow \infty}\frac{V''(x)}{ V'(x)^2}=0$$ satisfies the (MLS$(H_q)$) inequality if
and only if  $$\lim_{\vert x \vert\rightarrow \infty}\frac{V(x)}{\vert V' (x)\vert^q}<\infty$$
As a consequence, one obtains that non log-concave measures like for instance $$V(x)=\vert x \vert^p+\alpha (x)\cos(x)$$
such that $\alpha (x)=k\vert x \vert^{p-1-\delta}$ for $\vert x \vert >1$, $k$ a small constant, and $p\geq \frac{q}{q-1}$ and $\delta\in(0,1)$ satisfy the (MLS$(H_q)$) inequality. 

The purpose of this paper is to present examples where the (MLS$(H_q)$) inequality is satisfied for measures that go beyond the last mentioned example, as for instance  measures defined with  $$V(x)=d( x )^p+\alpha (x)\cos(x)$$
such that $\alpha (x)=kd( x )^{p-1}$ for $\vert x \vert >1$, $k\in (0,1)$ some small constant and $p\geq \frac{q}{q-1}$, where $d(x)=d(x_{0},x)$ denotes the distance of $x$ from a specific point $x_0$.   In order to obtain such examples, at first a criterion is going to be presented based on the U-bound inequalities introduced by Hebisch and Zegarlinski in \cite{H-Z}.  

In the aforementioned paper,
 the U-bound inequalities $$\mu(f^qd^s)\leq C\mu \vert \nabla f \vert^q +D \mu \vert f\vert^q$$where used to prove Log-Sobolev $q$ inequalities for $q\in (1,2]$, the spectral Gap inequality, as well as F-Sobolev inequalities.
In particular, for $\theta \in (1,2)$ the following U-bound inequality  \begin{align}\label{1stUbound}\mu_{\theta}(f^2d^{2(\theta-1)})\leq C\mu_{\theta} \vert \nabla f \vert^2 +D \mu_{\theta} \vert f\vert^2\end{align} was used to prove the F-Sobolev inequality
\begin{align}\label{F-Sobolev}\int f^2 \left\vert \log \frac{f^2}{\int f^2 d\mu_\theta} \right\vert^\varsigma d\mu_{\theta}\leq C\mu_{\theta} \vert \nabla f \vert^2 +D \mu _{\theta}\vert f\vert^2 \end{align}
for $\varsigma\in [0,1]$ and $\mu_\theta=\frac{e^{-d^{\theta}(x)}}{\int e^{-d^{\theta}(x)}dx}dx$.  

In this paper,  we apply the  methods developed in \cite{H-Z},  in order to obtain similar results for  the case of the Modified Logarithmic Sobolev inequality (MLS$(H_q)$).  More detailed, as shown in Theorem \ref{thm} the following U-bound inequality for $q>2$$$\mu(f^2d^{q(p-1)})\leq C\mu\left( \vert \nabla f \vert^q \vert f \vert^{2-q} \right)+D \mu \vert f\vert^2$$   is sufficient for the measure $\mu$ to satisfy an (MLS($H_q$)) inequality, for  $p\geq q'$, where $q'$ is the conjugate of $q$. As a consequence,  examples of non log-concave measures that satisfy a Modified Log-Sobolev inequality are presented in Corollary \ref{maincorollary}.
\section{A criterion on MLS$(H_q)$ inequalities.}
We consider general $n$-dimensional non compact metric spaces. We will denote $d$ the distance and  $\nabla$ the (sub)gradient for which we assume that $\frac{1}{\sigma}<\vert \nabla d \vert\leq 1$ for some $\sigma \in [1,\infty)$, and $\Delta d\leq K$ outside the unit ball $B=\{d(x)<1\}$ for some $K\in(0,+\infty)$. If  $d\lambda$ is the $n-$dimensional Lesbegue measure we assume that it satisfies the  Classical-Sobolev inequality (C-S)   
$$ \left(\int \vert f\vert ^{2+\epsilon}d\lambda\right)^\frac{2}{2+\epsilon}\leq\alpha \int\vert\nabla f  \vert^2d\lambda +\beta\int \vert f \vert ^2d\lambda  \   \   \   \   \   \   \   \   \   \       \text{(C-S)}$$ for positive constants $\alpha, \beta$, as well as    the local Poincar\'e inequality (see \cite{SC}),  that is, there exists a constant $c_R\in(0,\infty)$ such that for every ball $B_R$,\begin{align}\label{introSG1}\frac{1}{\vert B_R\vert}\int_{B_R}\left \vert f-\frac{1}{\vert B_R\vert}\int_{B_R}f\right \vert^2d\lambda\leq c_R\frac{1}{B_R} \int_{B_R} \vert \nabla f\vert^2d\lambda\end{align}
 \begin{Remark} The main assumptions for the $n$-dimensional non compact space with distance $d$ and (sub)gradient $\nabla$  are for instance satisfied in the case of the $\mathbb{R}^n,n\geq 1$ with $d$ being the Eucledian distance, as well as for the case of the Heisenberg group, with $d$ being the Carnot-Carath\'eodory distance.   

The local Poincar\'e inequality  for the Lebesgue measure (\ref{introSG1}) is a standard result (see for instance \cite{H}, \cite{H-Z}).  Concerning the  Classical Sobolev inequality for $n\geq 3$  one can look at \cite{D} and \cite{V-SC-C}.  For the  case of $n=1,2$ a stronger result  holds. The (C-S) in these two cases actually follows directly from the case $n=3$. To see this, for instance for $n=1$, if we consider a function $g:\mathbb{R}\rightarrow \mathbb{R}$    and then apply the (C-S) inequality for $n=3$ for the function $f(x_1,x_2,x_3)=\prod_{i=1}^3g(x_i)$,  we obtain $$\Vert g\Vert_{2+\epsilon}^6  \leq\alpha \left(\int\vert\nabla g   \vert^2d\lambda\right) \Vert g\Vert_{2}^4+\beta\Vert g\Vert_{2}^6$$
If we assume $\Vert g\Vert_2^2 =1$ we then have 
\begin{align}\label{proof(C-S)1}\Vert g\Vert_{2+\epsilon}^2  \leq\left(\alpha\int\vert\nabla g   \vert^2d\lambda+\beta\right)^\frac{1}{3}\end{align}
In the case where  $\int\vert\nabla g   \vert^2d\lambda\leq 1$,  (\ref{proof(C-S)1}) becomes 
\begin{align}\label{proof(C-S)2}\Vert g\Vert_{2+\epsilon}^2  \leq\left(\alpha+\beta\right)^\frac{1}{3}\end{align}
In the case where  $\int\vert\nabla g  \vert^2d\lambda> 1$, from (\ref{proof(C-S)1}) we obtain  \begin{align}\label{proof(C-S)3}\nonumber\Vert g\Vert_{2+\epsilon}^2 &  \leq\left(\alpha\int\vert\nabla g   \vert^2d\lambda+\beta\right)^\frac{1}{3}\leq\left(\alpha\int\vert\nabla g   \vert^2d\lambda\right)^\frac{1}{3}+\beta^\frac{1}{3}\\ &\leq \gamma \int\vert\nabla g   \vert^2d\lambda+\beta^\frac{1}{3} \end{align}
for positive constants $\gamma$ and $\beta$. If we combine together  inequalities  (\ref{proof(C-S)2}) and (\ref{proof(C-S)3}) we have that for any  $g:\mathbb{R}\rightarrow\mathbb{R}$ such that  $\Vert g\Vert_2^2 =1$ the following holds\ $$\Vert g\Vert_{2+\epsilon}^2\leq \gamma \int\vert\nabla g   \vert^2d\lambda+\left(\alpha+\beta\right)^\frac{1}{3}  $$ for positive constants $\alpha, \beta, \gamma$. The result follows if we replace  $g$ by $\frac{g}{\left(\int g^2d\lambda\right)^\frac{1}{2}}$.\end{Remark}
Furthermore,   for $d \lambda$ the Lesbegue measure,  we define  the probability measure \begin{align*}d\mu_{p}=\frac{e^{-d^{p}}}{Z_{p}}d\lambda
\end{align*}where $Z_{p}$ is the normalization constant. Since in this paper we are concerned with the subquadratic case we consider $1<p<2$. The main result is as follows.
\begin{theorem}\label{thm}For any $q>2$, let $d\mu=\frac{e^{-W}}{Z}d\mu_p$ for $p\geq\frac{q}{q-1}$ be a probability measure  defined with an a.e.  differentiable potential $W$ satisfying
$$\vert \nabla W\vert\leq \delta d^{p-1}+\gamma_\delta$$
with some small constant  $\delta \in (0,1)$ and $\gamma_\delta \in (0,\infty)$.   Then the following Modified  Log-Sobolev inequality (MLS($H_q$)) holds$$\mu\left\vert f\right\vert^2 log\frac{\left\vert f\right\vert^2}{\mu\left\vert f\right\vert^2}\leq c \ \int H_{q}\left(\frac{\vert \nabla f\vert}{f}\right)f^2d \mu$$for some positive constant $c$.\end{theorem}
As a direct consequence of the last theorem, we can obtain examples of measures that are not log-concave and actually go beyond the examples provided in \cite{B-R2}. The corollary below presents such a family of measures.
\begin{corollary}\label{maincorollary}Let $d\mu=\frac{e^{-V}}{\int e^{-V}d\lambda}d\lambda$ be a probability measure  defined with an a.e. differentiable potential $V$ satisfying
  $$V(x)=d( x )^p+\alpha (x)\cos(x)$$
such that $\alpha (x)=kd( x )^{p-1}$ for $\vert x \vert >1$,  where  $p\geq\frac{q}{q-1}$ and $k\in(0,1)$ a small constant.   Then the following Modified Log-Sobolev inequality holds$$\mu\left\vert f\right\vert^2 log\frac{\left\vert f\right\vert^2}{\mu\left\vert f\right\vert^2}\leq c  \int H_{q}\left(\frac{\vert \nabla f\vert }{ f}\right )f^2d\mu$$
 \end{corollary}
 As explained in the introduction, the $U-$bound inequalities introduced in \cite{H-Z} will play a crucial role in proving the MLS inequality. The proposition bellow provides a link between subquadratic measures and   $U-$bound inequalities. 
\begin{proposition}\label{Ubound1}  Let $d\mu=\frac{e^{-W}}{Z}d\mu_p$ be a probability measure defined with a differentiable potential $W$ satisfying
$$\vert \nabla W\vert\leq \delta d^{p-1}+\gamma_\delta$$
with  some small constant $\delta\in(0,1)$ and $\gamma_\delta \in (0,\infty)$. Then there exist constants $C',D'\in(0,+\infty)$ such that the following bound holds
\begin{align}\label{2ndUbound}\int \vert f \vert^2 d^{q(p-1)}d\mu\leq C\int H_q\left(\frac{ \vert \nabla f \vert}{ f  }\right)f^2 d\mu+D\int \vert f \vert^{2} d\mu\end{align}
 \end{proposition}
 
 \begin{proof} The starting point of the proof is the following U-bound inequality from  \cite{H-Z}.
\begin{theorem}\label{[H-Z]thm1}\textbf{(\cite{H-Z})} Assume that  $\frac{1}{\sigma}<\vert \nabla d \vert\leq 1$ for some $\sigma \in [1,\infty)$, and $\Delta d\leq K$ outside the unit ball $B=\{d(x)<1\}$ for some $K\in(0,+\infty)$. Let $d\mu=\frac{e^{-W}}{Z}d\mu_p$ be a probability measure defined with a differentiable potential $W$ satisfying
$$\vert \nabla W\vert\leq \delta d^{p-1}+\gamma_\delta$$
with some small constant $\delta\in(0,1)$ and $\gamma_\delta \in (0,\infty)$. Then there exist constants $C',D'\in(0,+\infty)$ such that the following bound holds
$$\int \vert f \vert d^{p-1}d\mu\leq C'\int \vert \nabla f \vert d\mu+D'\int \vert f \vert d\mu$$
 \end{theorem}
We will use the last theorem in order to obtain a $U-$bound inequality with a tighter left hand side. Let $d_1(x)=\max(1,d(x))$. From Theorem \ref{[H-Z]thm1}, by enlarging the constant  $D'$  we may assume that
\begin{equation}\label{new1U1}\int \vert f \vert d_{1}^{p-1}d\mu\leq C'\int \vert \nabla f \vert d\mu+D'\int \vert f \vert d\mu\end{equation}
If we choose $h=\vert f \vert^2d_1^{(q-1)(p-1)}$ we have
\begin{align}\label{new2U1}\int \vert f \vert^{2} d^{q(p-1)}d\mu&\leq\int \vert f \vert^{2} d_{1}^{q(p-1)}d\mu
\end{align}
Furthermore, from the inequality  (\ref{new1U1})  we can obtain
the following bound \begin{align}\label{U1.1}\int \vert f \vert^{2} d_{1}^{q(p-1)}d\mu=\int \vert h \vert d_{1}^{p-1}d\mu  \leq C'\int \vert \nabla h \vert d\mu+D'\int \vert h \vert d\mu\end{align}
For $R>1$ we have \begin{align}\nonumber \label{U1.6}\int h d\mu & = \int \vert f\vert^2d_1^{(p-1)(q-1)}d\mu\\ &=\int_{B_{R}} \vert f\vert^2d_1^{(p-1)(q-1)}d\mu+\int_{B_R^c} \vert f\vert^2d_1^{(p-1)(q-1)}d\mu \end{align}
where $B_R$ denotes a ball of radius $R$, i.e. $B_R=\{d(x)<R\}$. For the first term on the right hand side of (\ref{U1.6}) we have \begin{align}\label{U1.7}
\int_{B_R} \vert f\vert^2d_1^{(p-1)(q-1)}d\mu\leq R^{(p-1)(q-1)}\int \vert f\vert^2d\mu\end{align}
While for the second term on the right hand side of (\ref{U1.6})  we compute
\begin{align}\label{U1.8}\int_{B_R^c}\vert f\vert^2d_1^{(p-1)(q-1)}d\mu\leq\frac{1}{R^{p-1}} \int \vert f\vert^2d_1^{q(p-1)}d\mu
\end{align}
From (\ref{U1.6})-(\ref{U1.8}) we finally obtain
\begin{align}\label{U1.9}\int h d\mu =R^{(p-1)(q-1)}\int \vert f\vert^2d\mu+\frac{1}{R^{p-1}} \int \vert f\vert^2d_1^{q(p-1)}d\mu
\end{align}
for a constant $R>1$. Furthermore, we have 
 \begin{align}\label{U1.2} \vert \nabla h\vert =2\vert \nabla f \vert \vert f \vert d_1^{(p-1)(q-1)}+(q-1)(p-1)\vert \nabla d_1 \vert \vert f\vert^2 d_1^{(p-1)(q-1)-1} \end{align}
We can compute 
\begin{align}\label{U1.3}\nonumber\int 2\vert \nabla f \vert  \vert f \vert &d_1^{(p-1)(q-1)}d\mu \\ & \nonumber=\int 2\vert \nabla f \vert\vert f\vert^{\frac{2}{q}-1} \vert f \vert^{2\frac{(q-1)}{q}} d_1^{(p-1)(q-1)}d\mu\\  & \nonumber\leq2\left( \int \vert \nabla f \vert^{q}\vert f\vert^{2-q} d\mu \right)^{\frac{1}{q}}\left( \int \vert f\vert^{2}  d_1^{\frac{q(p-1)(q-1)}{q-1}}d\mu \right)^{\frac{q-1}{q}} \\  &    \leq\alpha^2\int \vert \nabla f \vert^{q}\vert f\vert^{2-q} d\mu +\frac{1}{\alpha^2}\int \vert f\vert^{2}  d_1^{q(p-1)}d\mu
\end{align}
for a constant $\alpha>0$. We also have 
\begin{align}\label{U1.4}\nonumber\int\vert \nabla d_1 \vert \vert f\vert^2 d_1^{(p-1)(q-1)-1}d\mu  & \leq\int f^2 d_1^{(p-1)(q-1)}d\mu \\  &  \leq R^{(p-1)(q-1)}\int f^2d\mu+\frac{1}{R^{p-1}} \int_{} \vert f\vert^2d_1^{q(p-1)}d\mu
\end{align}where above we used the bound from (\ref{U1.9}). From (\ref{U1.2}), (\ref{U1.3}) and (\ref{U1.4}) we obtain 
\begin{align}\nonumber \label{U1.5}\int \vert \nabla h \vert d\mu\leq&\alpha^2\int \vert \nabla f \vert^{q}\vert f\vert^{2-q} d\mu+(q-1)(p-1)R^{(p-1)(q-1)}\int \vert f\vert^2d\mu \\ &+\left(\frac{1}{\alpha^2}+\frac{(q-1)(p-1)}{R^{p-1}}\right) \int_{} \vert f\vert^2d_1^{q(p-1)}d\mu 
\end{align} 
If we plug (\ref{U1.5}) and (\ref{U1.9}) in (\ref{U1.1}) we finally obtain
\begin{align}\label{U1.10}\nonumber\int \vert f \vert^2 d_{1}^{q(p-1)}d\mu\leq  & C'\alpha^2\int \vert \nabla f \vert^q\vert f\vert^{2-q} d\mu \\ & \nonumber+\left(D'+C'(q-1)(p-1)\right) R^{(p-1)(q-1)}\int \vert f\vert^2d\mu \\ & +\left(\frac{C'}{\alpha^2}+\frac{D'+C'(q-1)(p-1)}{R^{p-1}}\right) \int \vert f\vert^2d_1^{q(p-1)}d\mu \end{align}
If we choose $\alpha$ and $R$ large enough so that $\frac{C'}{\alpha^2}+\frac{D'+C'(q-1)(p-1)}{R^{p-1}}<1$ we obtain 
\begin{align}\label{newlastU1} \int \vert f \vert^{2} d_1^{q(p-1)}d\mu\leq C\int \vert \nabla f \vert^{q}\vert f\vert^{2-q} d\mu+\breve  D\int  f^2d\mu  \end{align}
for constants $$ C=\frac{C'\alpha}{1-\left(\frac{C'}{\alpha^2}+\frac{D'+C'(q-1)(p-1)}{R^{p-1}}\right) }$$ and $$\breve D=\frac{\left(D'+C'(q-1)(p-1)\right)R^{(p-1)(q-1)}}{1-\left(\frac{C'}{\alpha^2}+\frac{D'+C'(q-1)(p-1)}{R^{p-1}}\right) }$$ The proof of the proposition follows from (\ref{new2U1}) and  (\ref{newlastU1}) for constant $D=C+\breve D$, 
since $\vert \nabla f \vert^{q}\vert f\vert^{2-q} \leq f^2$ when $\vert \nabla f\vert\leq \vert f \vert$.   \end{proof}

 ~

 If we compare the U-bound inequality (\ref{2ndUbound}) of Proposition \ref{Ubound1} with the U-bound inequality (\ref{1stUbound}) used in \cite{H-Z} to show the F-Sobolev inequality (\ref{F-Sobolev}), one notices that the left hand side of (\ref{2ndUbound}) is stronger, while the right hand side is relaxed from the full gradient to the weaker modification related with the Modified Logarithmic Sobolev inequality (MLS($H_q$)). In the next proposition we present the link between  the $U$-bound   inequality of the last proposition and the Defective Modified Log-Sobolev inequality. 
 \begin{proposition}\label{Ubound2} Suppose that the measure $$d\mu=\frac{e^{-U}d\lambda}{\int e^{-U}d\lambda\ }$$ 
 where $d\lambda$ the $n-$dimensional  Lebesgue measure and $U\geq 0$, satisfies the following U-bound inequality \begin{equation}\label{Uboundprop1.2p5}\mu\vert f\vert ^2\left(\vert \nabla U\vert^2+U\right)\leq\hat C\int H_{q}\left(\frac{\vert \nabla f\vert }{ f}\right )f^2d\mu+\hat D\mu\vert f \vert ^2\end{equation}
for some  positive constants $\hat C$ and   $\hat D$ both independent of $f$. Then the following Defective Modified Log-Sobolev inequality holds$$\mu\left\vert f\right\vert^2 log\frac{\left\vert f\right\vert^2}{\mu\left\vert f\right\vert^2}\leq C \int H_{q}\left(\frac{\vert \nabla f\vert }{ f}\right )f^2 d\mu +D\mu f^2$$
\end{proposition}
\begin{proof}
We follow closely the work in \cite{H-Z}  for the Log-Sobolev $q$ inequality (see also   \cite{I-P}). Without loss of generality we can assume that $f\geq 0$ and we set $\rho=\frac{e^{-U}}{\int e^{-U}d\lambda }$ and $g=f\rho ^\frac{1}{2}$ We also assume that $$\int g^2d\lambda=\mu f^2=1$$
Then we can write 
$$\int (g^2\log g^2)d\lambda=\frac{2}{\epsilon}\int g^2(\log 2^\epsilon)d\lambda_{}\leq \frac{2+\epsilon}{2}\frac{2}{\epsilon}\log\left(\int g^{2+\epsilon}d\lambda\right)^\frac{2}{2+\epsilon}$$
where above we used the Jensen's inequality. If we use now the Classical-Sobolev inequality (C-S) for the   Lebesgue measure $d\lambda$  
$$ \left(\int \vert f\vert ^{2+\epsilon}d\lambda\right)^\frac{2}{2+\epsilon}\leq\alpha \int\vert\nabla f  \vert^2d\lambda +\beta\int \vert f \vert ^2d\lambda  \   \   \   \   \   \   \   \   \   \       \text{(C-S)}$$ for positive constants $\alpha, \beta$, we will get
\begin{align}\nonumber\int (g^2\log g^2)d\lambda&\leq\frac{2+\epsilon}{\epsilon}\log\left(\alpha \int\vert\nabla g \vert^2d\lambda +\beta\int \vert g \vert ^2d\lambda\right) \\ &\label{eq3prop1.2p5}\leq\frac{(2+\epsilon)\alpha}{\epsilon} \int\vert\nabla g \vert^2d\lambda+\frac{(2+\epsilon)\beta}{\epsilon}\int \vert g\vert^2d\lambda\end{align}
where in the last inequality we used that $\log x\leq x$ for $x>0$.  For the first term on the right hand side of (\ref{eq3prop1.2p5}) we have
\begin{align}\nonumber\label{eq4prop1.2p5}\int \vert \nabla g \vert^2d\lambda&=\int \vert \nabla (f\rho^\frac{1}{2})  \vert^2d\lambda\\  &\leq2^{q-1}\mu\vert \nabla f  \vert^2+2\int \vert f\nabla (\rho^\frac{1}{2})  \vert^2d\lambda\end{align}
We have \begin{align*}\int \vert f\nabla (\rho^\frac{1}{2})  \vert^2d\lambda&=\int \vert \rho^\frac{1}{2}\rho^\frac{-1}{2}f\nabla (\rho^\frac{1}{2})  \vert^2d\lambda=\mu f^{2} \vert\rho^\frac{-1}{2} \nabla (\rho^\frac{1}{2})  \vert^2\\ &=\frac{1}{4}\mu f^{2} \vert\nabla U\vert^2\end{align*}
If we plug the last equality in (\ref{eq4prop1.2p5}), we obtain 
\begin{equation}\label{eq5prop1.2p5} \int \vert \nabla g \vert^2d\lambda\leq2\mu\vert \nabla f  \vert^2+\frac{1}{2}\mu f^{2} \vert\nabla U\vert^2
\end{equation}
If we combine inequalities  (\ref{eq3prop1.2p5}) and (\ref{eq5prop1.2p5}),  we get
\begin{align}\nonumber\int (g^2\log g^2)d\lambda\leq&\frac{2(2+\epsilon)\alpha}{\epsilon} \mu\vert \nabla f  \vert^2+\frac{(2+\epsilon)\beta}{\epsilon}\mu f^2\\ \label{eq6prop1.2p5}&+\frac{2(2+\epsilon)\alpha}{\epsilon4}\mu f^{2} \vert\nabla U\vert^2\end{align}
For the left hand side of (\ref{eq6prop1.2p5}), since $U\geq 0$,  we have 
\begin{align}\int (g^2\log g^2)d\lambda\nonumber=&\int \left(\frac{e^{-U}}{\int e^{-U}d\lambda\ }f^2\log \frac{e^{-U}}{\int e^{-U}dX }f^2\right)d\lambda\\=&\mu (f^2\log f^2)+\nonumber\mu \left(f^2\log \frac{e^{-U}}{\int e^{-U}d\lambda }\right)
\\=&\label{eq7prop1.2p5}\mu (f^2\log f^2)-\mu(f^2U) \nonumber-\mu \left(f^2\log \int e^{-U}d\lambda\right)\\\geq&\mu (f^2\log f^2)-\mu(f^2 U)\end{align}
  If we combine (\ref{eq6prop1.2p5}) and (\ref{eq7prop1.2p5}), we obtain
\begin{align}\label{eq1.8prop1.2p5}\mu(f^2\log f^2)\leq&\hat\alpha\mu\vert \nabla f  \vert^2+\hat \gamma\mu f^2+\hat \beta\mu f^{2} \left(\vert\nabla U\vert^2+U\right)\end{align}
where $\hat \alpha=\frac{2(2+\epsilon)\alpha}{\epsilon}$, $\hat \beta=\max\{\frac{2(2+\epsilon)\alpha}{\epsilon 4},1\}$ and $\hat \gamma =\frac{(2+\epsilon)\beta}{\epsilon}$. If we use the  U-bound  (\ref{Uboundprop1.2p5}),  the inequality (\ref{eq1.8prop1.2p5}) gives
$$\mu(f^2\log f^2)\leq(\hat\alpha+\hat \beta\hat C)\int H_{q}\left(\frac{\vert \nabla f\vert }{ f}\right )f^2d\mu+\hat \beta\hat D\mu f^2+\hat \gamma\mu f^2$$
If we replace $f$ with $\frac{f}{\mu f }$ which has mean equal to one we obtain the result. 
\end{proof}

~

We can now present the proof of Theorem \ref{thm}.

\noindent 
\textbf{\textit{Proof of Theorem \ref{thm}}.} Since $p\geq \frac{q}{q-1}$  and $\vert \nabla W\vert\leq \delta d^{p-1}+\gamma_\delta$ for a small constant $\delta\in(0,1)$ and $\gamma_\delta \in (0,\infty)$ we derive that for $U=d^p+W$ there exists a positive constant $\tilde C$ such that $$\int f^2\left(\vert \nabla U\vert^2+U\right)d\mu\leq  \tilde C \int f^2d^{q(p-1)}d\mu+\tilde C\int f^2 d\mu$$since $q>2$ and $p\geq \frac{q}{q-1}$. If we apply the U-bound inequality of Proposition \ref{Ubound1} to bound the first term on the right hand side of the last inequality, we then get 
\begin{align*}\int f^2\left(\vert \nabla U\vert^2+U\right)d\mu\leq\tilde  CC\int H_q\left(\frac{ \vert \nabla f \vert}{ f  }\right)f^2 d\mu+\tilde C(D+1)\int \vert f \vert^{2} d\mu
\end{align*} Since we have obtained a U-bound inequality like the one required  in  hypothesis  (\ref{Uboundprop1.2p5})  we can   apply Proposition \ref{Ubound2}. This will lead to the  Defective Modified Log-Sobolev inequality   \begin{align}\label{proof1}\int f^2 \log \frac{f^2}{\int f^2 d\mu} d\mu\leq C \int f^2H_{q}\left(\frac{\vert \nabla f \vert}{f}\right )d\mu+D\int f^2 d\mu\end{align} for positive constants $C$ and $D$.  In order to finish the  proof of the theorem it remains to pass from the Defective Modified Log-Sobolev inequality (\ref{proof1}) to the Modified Log-Sobolev inequality (MLS($H_q$)). To do this we will first need the Poincar\'e inequality for the measure $d\mu$. This is provided from the  following  theorem, whose proof can be found in  \cite{H-Z}. 
\begin{theorem}\label{thmSG[H-Z]}\textbf{(\cite{H-Z})} Suppose $1\leq q<\infty$ and a measure $\lambda$ satisfies the $q-$Poincar\'e inequality for every ball $B_R$, that is there exists a constant $c_R\in(0,\infty)$ such that \begin{align}\label{SG1}\frac{1}{\vert B_R\vert}\int_{B_R}\left \vert f-\frac{1}{\vert B_R\vert}\int_{B_R}f\right \vert^qd\lambda\leq c_R\frac{1}{B_R} \int_{B_R}\vert \nabla f\vert^qd\lambda\end{align}
Let $\mu$ be a probability measure which is absolutely continuous with respect to the measure $\lambda$ and such that 
\begin{align}\label{SG2}\int \vert f\vert ^q\eta d\mu\leq C\int \vert \nabla f\vert ^q+D\int \vert f \vert^qd\mu
\end{align}
with some nonnegative function $\eta$ and some constants $C,D\in (0,\infty)$ independent of a function $f$.  If for any $L\in(0,\infty)$ there is a constant $A_L$ such that 
\begin{align*}\frac{1}{A_L}\leq \frac{d\mu}{d\lambda}\leq A_L
\end{align*}
on the set $\{\eta<L\}$ and, for some $R\in(0,\infty)$ (depending on $L$), we have $\{\eta<L\}\subset B_R$, then $\mu$ satisfies the q-Poincar\'e inequality
\begin{align*}\mu \left \vert f-\mu f \right \vert ^q\leq c' \mu \left \vert \nabla f \right \vert^q
\end{align*}
\end{theorem}  For 
$\lambda$ being the Lebesgue measure, (\ref{SG1}) is true as initially  assumed in (\ref{introSG1}). Furthermore,  (\ref{SG2}) was shown in Proposition  \ref{Ubound1}  for $\eta=d^{q(p-1)}$. Thus,   if we apply the above theorem for $q=2$   we obtain  that the measure $\mu$ also satisfies the following Spectral Gap inequality 
\begin{align}\label{proof2}\mu \left \vert f-\mu f \right \vert ^2\leq c' \mu \left \vert \nabla f \right \vert ^2
\end{align}
for some positive constant $c'$. If we combine the Defective Modified Logarithmic Sobolev inequality   (\ref{proof1}), together with the Spectral Gap inequality (\ref{proof2}) the theorem follows according to the following theorem (see \cite{B-K}) 
\begin{theorem}\label{thm[B-K]}\textbf{(\cite{B-K})} Let $H$ be an even function on $\mathbb{R}$, which is increasing on $\mathbb{R}^+$ and satisfies $H(0)=0$ and $H(x)\geq cx^2$. Assume that there exists $q\geq 2$ such that $x\rightarrow\frac{H(x)}{x^q}$ is non-increasing on $(0,+\infty)$.

 Assume that a probability measure $\mu$ satisfies a Defective Modified Log-Sobolev inequality:$$\int f^2 \log\frac{f^2}{\int f^2 d\mu}d\mu\leq \int f^2H\left(\frac{\vert \nabla f \vert}{f}\right )d\mu+D\int f^2 d\mu$$
If $\mu$ also satisfies a Poincar\'e inequality, then there exists a constant $C$ such that for every $f$,
 $$\int f^2 \log\frac{f^2}{\int f^2 d\mu}d\mu\leq C \int f^2H\left(\frac{(\vert \nabla f \vert}{f}\right )d\mu$$\end{theorem}
\qed
 
\section{Conclusion.}
In this paper a criterion on the Modified Log-Sobolev inequality was presented with the use of U-bounds.  In particular, in Proposition \ref{Ubound2}  we saw that the following U-bound inequality  \begin{equation*}\mu\vert f\vert ^2\left(\vert \nabla U\vert^2+U\right)\leq\hat C\mu \left(\ H_{q}\left(\frac{\vert \nabla f\vert }{ f}\right )f^2\right)+\hat D\mu\vert f \vert ^2\end{equation*}together with the Spectral Gap are sufficient for the measure $d\mu=\frac{e^{-U}d\lambda}{\int e^{-U}d\lambda}d\lambda$ to satisfy a MLS($H_q$) inequality. Concerning the converse problem, if we follow the work in \cite{H-Z} for the similar problem concerning the stronger Log-Sobolev $q$ inequalities we can show the following result.
\begin{theorem} Suppose that the measure $d\mu=\frac{e^{-U}d\lambda}{\int e^{-U}d\lambda}d\lambda$ satisfies an MLS($H_q$) inequality and that $U$ is such that
$$\vert \nabla  U \vert^q\leq aU+b$$
for positive constants $a$ and $b$. Then the following U-bound is true
$$\int\vert f\vert ^2 Ud\mu\leq C\int H_{q}\left(\frac{\vert \nabla f\vert }{ f}\right )f^2d\mu+ D\int\vert f \vert ^2\mu$$
\end{theorem}
This demonstrates that as in the case of Log-Sobolev $q$  and F-Sobolev inequalities, the Modified Log-Sobolev MLS($H_q$) inequalities are also equivalent to  U-bound inequalities.


\end{document}